\documentclass[11pt]{amsart}
\usepackage{amssymb}
\usepackage[cmtip,matrix,arrow]{xy}

\newtheorem{theorem}{Theorem}[section]

\newtheorem{proposition}[theorem]{Proposition}

\newtheorem{lemma}[theorem]{Lemma}

\theoremstyle{remark}
\newtheorem{remark}[theorem]{Remark}

\newcommand\be{\begin{equation}\label}
\newcommand\ee{\end{equation}}

\newcommand\lie[1]{\mathfrak{#1}}

\newcommand{\g}{\lie{g}}

\newcommand\dirac{/\kern-1.2ex\partial} % Dirac operator
\newcommand\qu{/\kern-.7ex/} % Categorical quotients
\newcommand{\hra}{\hookrightarrow}

\newcommand{\p}{\partial}
\renewcommand{\l}{\langle}
\renewcommand{\r}{\rangle}

\newcommand{\ti}{\tilde}

\newcommand\beqn{\begin{equation}}      
\newcommand\eeqn{\end{equation}}

\newcommand{\beq}{\begin{eqnarray*}}
\newcommand{\eeq}{\end{eqnarray*}}

\begin{document}

\title[ ]{Poisson geometry and \\
the Kashiwara-Vergne conjecture}

\author{A. Alekseev}
\address{University of Geneva, Section of Mathematics,
2-4 rue du Li\`evre, 1211 Gen\`eve 24, Switzerland}
\email{alekseev@math.unige.ch}

\author{E. Meinrenken}
\address{University of Toronto, Department of Mathematics,
100 St George Street, Toronto, Ontario M5S3G3, Canada }
\email{mein@math.toronto.edu}

%\date{\today}

\begin{abstract}
We give a Poisson-geometric proof of the Kashiwara-Vergne conjecture for 
quadratic Lie algebras, based on the equivariant Moser trick.
\end{abstract}


\maketitle

\vskip 0.3cm

\section{Introduction}
In 1977 Duflo \cite{du:op} established the local solvability of
bi-invariant differential operators on arbitrary finite-dimensional
Lie groups $G$.  In a subsequent paper, Kashiwara-Vergne
\cite{ka:ca} conjectured a property of the Campbell-Hausdorff series
that would imply Duflo's result as well as a more general
statement on convolution of (germs of) invariant distributions.
They proved this property (cf. Section
\ref{sec:interpret} below) for solvable Lie algebras.

Except for case $\g={\mathfrak{sl}}(2,{\mathbb R})$ proved by Rouvi\`{e}re
\cite{ro:de}, there had not been much progress on the conjecture for
more than twenty years until Vergne \cite{ve:ce} settled the case of
quadratic Lie algebras.  Shortly after, results for the general case
were obtained in \cite{an:co,an:de,to:su} using Kontsevich's approach
to deformation quantization.

In this note we will give a new proof of the Kashiwara-Vergne
conjecture for quadratic Lie algebras, using ideas from Poisson
geometry. We will introduce a family of Poisson structures $P_t$ on a
neighborhood of the origin in $\g\times\g$, with the property that the
diagonal action of $G$ is Hamiltonian, with moment map
$\Phi_t(X,Y)=\frac{1}{t}\log(e^{tX}e^{tY})$.  Using Moser's argument we
will define a time dependent vector field $v_t$ interpolating the
Poisson structures and moment maps, and in particular relating
addition $\Phi_0=X+Y$ and multiplication $\Phi_1=\log(e^X e^Y)$. We
will see that Moser's vector field directly gives a solution to the
Kashiwara-Vergne problem.

{\bf Acknowledgements.}
We would like to thank Michele Vergne for explaining to us the ideas
behind \cite{ka:ca} and \cite{ve:ce}, and helpful discussions.

\section{Geometric formulation of the Kashiwara-Vergne property}
\label{sec:interpret}
Let $\g$ be a finite dimensional Lie algebra, and $G$ the
corresponding simply connected Lie group. Choose an open neighborhood
$U\subset \g$ of $0$, such that the exponential map $\exp:\,\g\to
G,\,X\mapsto e^X$ is a diffeomorphism over $U$. Denote by
$\log:\,\exp(U)\to \g$ the inverse.  Let $V\subset \g\times \g$ be a
convex open neighborhood of the origin, with the property that $e^X
e^Y\in\exp(U)$ for $(X,Y)\in V$.  Given $F:\,V\to \g$ denote by
$F\delta_X$ the vector field $f\mapsto \frac{d}{d s}|_{s=0}f(X+s
F(X,Y),Y),\ \ f\in C^\infty(V),$ and by $\delta_X f$ the
${\rm End}(\g)$-valued function, $\ a\mapsto \frac{d}{d s}|_{s=0}f(X+sa,Y),\ \
a\in\g$. Similarly, define $F\delta_Y$ and $\delta_Y f$.
\vskip.2in

\noindent{\bf Conjecture (Kashiwara-Vergne)} There exist $\g$-valued analytic
functions $A,B$ on a neighborhood of the origin of
$\g\times\g$, with $A(0,0)=B(0,0)=0$, such that
\begin{equation}
\label{eq:1}
\log(e^Y\,e^X)-X-Y=
(1-e^{-{\rm ad}_X})A(X,Y)+(e^{{\rm ad}_Y}-1)B(X,Y)
\end{equation}
and
\begin{equation}
\label{eq:2}
{\rm tr}({\rm ad}_X \delta_X(A)+ {\rm ad}_Y \delta_Y(B))
=-\frac{1}{2} {\rm tr}\left(\frac{{\rm ad}_X}{e^{{\rm ad}_X}-1}+\frac{{\rm ad}_Y}{e^{{\rm ad}_Y}-1}-
\frac{{\rm ad}_Z}{e^{{\rm ad}_Z}-1}-1\right)
\end{equation}
where $Z=\log(e^X e^Y)$.

Kashiwara-Vergne arrived at these conditions from simple and natural
geometric considerations. They considered a deformation of the Lie
bracket $[\cdot,\cdot]_t=t[\cdot,\cdot]$ and examined the limit $t\to 0$
where the Lie algebra becomes Abelian. Retracing their argument, one
has the following equivalent formulation of the Conjecture.
Consider the vector field $v$ and the $\g$-valued function $\Phi$
on $V$ given by
\begin{equation}
\label{eq:v}
v=({\rm ad}_X A)\delta_X+({\rm ad}_Y B)\delta_Y,\ \
\Phi(X,Y)=\log(e^X e^Y).
\end{equation}
Let $m_t:\,V\to V$ denote multiplication by $t\in [0,1]$, and
rescale $v_t=\frac{1}{t}m_t^*v$ and $\Phi_t=\frac{1}{t}m_t^*\Phi$. (Note that
this is well-defined even for $t=0$, and that $\Phi_0(X,Y)=X+Y$.)
Let $J(X)=\det\left(\frac{1-e^{-{\rm ad}_X}}{{\rm ad}_X}\right)$ be the
Jacobian of the exponential map, and denote
by $\kappa_t\in C^\infty(V)$ the combination,
\begin{equation}
\label{eq:3a} \kappa_t(X,Y)=\frac{J^{1/2}(tX)J^{1/2}(tY)}{J^{1/2}(t\Phi_t)}.
\end{equation}
Let $\Gamma$ denote a translation invariant volume form on
$\g\times\g$.

\begin{proposition}[Kashiwara-Vergne]
Suppose $\g$ is unimodular. Then \eqref{eq:1}, \eqref{eq:2} are
respectively equivalent to the statements that the
vector field $v_t$ intertwines the maps $\Phi_t$ and the volume forms
$\kappa_t\Gamma$:
\begin{equation} 
\label{eq:1'} (\frac{\p}{\p t}+v_t)\Phi_t=0,\ \
(\frac{\p}{\p t}+v_t)(\kappa_t\Gamma)=0.
\end{equation}
\end{proposition}

For the proof, see \cite{ka:ca}, p.255--257. The unimodularity
assumption enters the reformulation of \eqref{eq:2}, since
the left hand side is the divergence of $v$ in this case.

\section{Poisson geometry}
Our approach to the Kashiwara-Vergne conjecture requires some elementary
concepts from Poisson geometry, which we briefly review in this Section.
%See e.g. the book \cite{ca:ge} for more detailed information.

\subsection{Basic definitions}
A Poisson manifold is a manifold $M$, equipped with a bi-vector field
$P\in C^\infty(M,\wedge^2 TM)$ such that the Schouten bracket $[P,P]$
vanishes. Let $P^\sharp:\,T^*M\to TM$ be the bundle map defined by
$P$, i.e. $P^\sharp(a)=P(a,\cdot)$ for all covectors $a$. The Poisson
structure gives rise to a generalized foliation of $M$, such that the
tangent bundle to any leaf $L\subset M$ equals
$P^\sharp(T^*M)|_L$. Each leaf carries a symplectic form $\omega_L$ given
as the inverse of $P|_L$.  Let $\Lambda_L$ denote the symplectic
volume form on $L$.

For any function $H\in C^\infty(M)$ one defines the Hamiltonian vector
field $v_H=-P^\sharp({\rm d} H)$. An action of a Lie group $G$ on $M$ is
called Hamiltonian if there exists an equivariant
{\em moment map} $\Phi\in
C^\infty(M,\g^*)$ with $ X_M=-v_{\l\Phi,X\r}$,
where $X_M=\frac{d}{ d t}|_{t=0}\exp(-t X)$ is the vector field generated
by $X\in\g$.  Equivalently, the restriction of $\Phi$ to any
symplectic leaf is a moment map in the sense of symplectic geometry.
For any Lie group $G$, the dual of the Lie algebra $\g^*$
carries a unique {\em Kirillov Poisson structure} such that
the identity map is a moment map for the co-adjoint action. Its symplectic
leafs are the co-adjoint orbits.
For any volume form $\Gamma$ on a Poisson manifold $(M,P)$, one defines
the {\em modular vector field} by $w_\Gamma(H)=-{\rm div}_\Gamma(v_H)$.
The modular vector field for $M=\g^*$ is given by modular character:
That is, $w=\sum_a {\rm tr}({\rm ad}(e_a))\frac{\p}{\p \mu_a}$ where $e_a$ is a
basis of $\g$ and $\mu_a$ the dual coordinates on $\g^*$. In particular,
$w=0$ for unimodular Lie algebras.

\begin{lemma}\label{lem:div}
Let $M$ be a Poisson manifold and $\Gamma$ be a volume form on $M$.
Suppose that the  modular vector field vanishes.
Then,  for any vector field $v$ of the form  $v=\sum_i F_i v_{H_i}$,
with smooth functions $F_i,H_i$, one has
$$  {\rm div}_\Gamma(v)|_L={\rm div}_{\Lambda_L}(v|_L).$$
for all symplectic leafs $L$.
\end{lemma}
\begin{proof}
We have ${\rm div}_{\Lambda_L}(v_{H_i}|_L)=0$ since Hamiltonian vector
fields preserve the symplectic form, and we have
${\rm div}_\Gamma(v_{H_i})|_L=0$ since the modular vector field vanishes.
Hence, both sides equal  the restriction of $\sum_i v_{H_i}(F_i)$ to $L$.
\end{proof}

\subsection{Gauge transformations}
There is an easy way of constructing new Hamiltonian Poisson
structures out of given ones referred to as `gauge
transformations' in Poisson geometry. It will be convenient to state this
using the language of equivariant de Rham theory. Let $M$ be a
$G$-manifold, where $G$ is a connected Lie group. One defines a
complex of equivariant differential forms $(\Omega_G^\bullet(M),{\rm d}_G)$,
where $\Omega_G(M)$ is the space of polynomial maps $\alpha:\,\g\to \Omega(M)$
satisfying the equivariance condition,
$\alpha([Y,X])+L_{Y_M}\alpha(X)=0$.  The equivariant differential is
defined as $({\rm d}_G\alpha)(X)= {\rm d}\alpha(X)-\iota_{X_M}\alpha(X)$. One
introduces a grading on $\Omega_G(M)$ by declaring that
${\deg}(\alpha)= 2k+l$ if $X\mapsto \alpha(X)$ is a homogeneous
polynomial of degree $k$ with values in $\Omega^l(M)$. In particular,
equivariant 2-forms are sums $\sigma_G(\xi)=\sigma+\l\Psi,\xi\r$ where
$\sigma$ is an invariant 2-form and $\Psi$ an equivariant map into
$\g^*$.  Note that these definitions extend to {\em local} $G$-actions
(equivalently, Lie algebra actions).

Suppose $(M,P_0,\Phi_0)$ is a Hamiltonian Poisson $G$-manifold, and
$\sigma_G=\sigma+\Psi\in\Omega^2_G(M)$ an equivariant cocycle.  Let
$\sigma^\flat:\,TM\to T^*M$ be the bundle map defined by $\sigma$.  Assume
that $\det(1+\sigma^\flat\circ P_0^\sharp)>0$ everywhere.  Then there is
a well-defined bivector field $P_1$ on $M$ such that
\begin{equation}
\label{eq:change}
 P_1^\sharp=P_0^\sharp\circ (1+\sigma^\flat\circ
 P_0^\sharp)^{-1}.\end{equation}
Let $\Phi_1=\Phi_0-\Psi$. Then $(M,P_1,\Phi_1)$ is again a Hamiltonian
Poisson $G$-manifold, with the same symplectic leafs as $P_0$. For any
leaf $L$, the equivariant symplectic forms $(\omega_j)_G=\omega_j-\Phi_j|_L$
are related by $(\omega_1)_G=(\omega_0)_G+\iota_L^* \sigma_G$, and the volume
forms by $\Lambda_{L,1}=\lambda|_L\Lambda_{L,0}$ where
$\lambda={\det}^{1/2}(1+\sigma^\flat\circ P_0^\sharp)$.

\subsection{Moser's trick for Poisson manifolds}
Let $M$ be a $G$-manifold, with a family of Hamiltonian
Poisson structures $(P_t,\Phi_t)$
depending smoothly on $t$, where the symplectic foliation is
independent of $t$. Suppose there exists a family of invariant 1-forms
$\alpha_t\in\Omega^1(M)^G$, depending smoothly on $t$, such that for every
symplectic leaf $\iota_L:\,L\hra M$,
$$ d_G( \iota_L^*\alpha_t)=\frac{d (\omega_t)_G}{d t}.$$
We will call $v_t=P_t^\sharp(\alpha_t)$ the {\em Moser vector field}.
Restricting to symplectic leafs, the usual
equivariant Moser's trick
from symplectic geometry (cf. \cite[Theorem 7.3]{ca:le} or \cite[Lemma 3.4]{al:li}) shows that $v_t$ intertwines the Poisson structures and moment maps:
$$(\frac{\p}{\p t}+v_t)\Phi_t=0,\ \ (\frac{\p}{\p t}+v_t)P_t=0 .$$
Furthermore, if $\Gamma_0$ is a volume form on $M$ with vanishing 
modular vector field with respect to $P_0$, a straightforward 
calculation, using Lemma \ref{lem:div}, shows that 
$$
(\frac{\p}{\p t}+v_t)\Gamma_t=0,
$$
where $\Gamma_t=\lambda_t \Gamma_0$ with $\lambda_t=
{\det}^{1/2}(1+\sigma_t^\flat\circ P_0^\sharp)$.

\section{Proof of the Kashiwara-Vergne conjecture for quadratic $\g$}

\subsection{Deformation of the Kirillov Poisson structure}\label{subsec:sig}
Suppose $\g$ is a quadratic Lie algebra. That is, $\g$ comes equipped
with a non-degenerate invariant symmetric bilinear form $\cdot$, used
to identify $\g^*\cong\g$. Quadratic Lie algebras are unimodular.  Let
$P_0$ denote the Poisson structure on $\g\times\g$ given as a product
of Kirillov Poisson structures. The moment map for the diagonal
$G$-action is the sum $\Phi_0(X,Y)=X+Y.$ Our goal is to find an
equivariantly closed equivariant 2-form
$\sigma_G=\sigma+\Psi\in\Omega^2_G(V)$ such that $\Psi=\Phi_0-\Phi_1$, where
$\Phi_1(X,Y)=\log(e^X e^Y).$ The construction is motivated by ideas
from the theory of group-valued moment maps, introduced in
\cite{al:mom}.
Let $\theta^L,\theta^R\in\Omega^1(G,\g)$ denote the left/right
invariant Maurer-Cartan on $G$. Let $G$ act on itself by
conjugation, and let
$\eta_G\in\Omega^3_G(G)$ be the equivariant Cartan 3-form
$$ \eta_G(X)=\frac{1}{12} \theta^L\cdot [\theta^L,\theta^L]
-\frac{1}{2} (\theta^L+\theta^R)\cdot X.
$$
It is well-known that ${\rm d}_G\eta_G=0$.
One easily calculates the pull-back of this form under group multiplication
${\rm Mult}:\,G\times G\to G$:
\begin{equation}
\label{eq:Mult} {\rm Mult}^*\eta_G=\eta_G^1+\eta_G^2-\frac{1}{2} {\rm d}_G
(\theta^{L,1}\cdot\theta^{R,2})\end{equation}
where the superscripts $1,2$ denote pull-back to the respective
$G$-factor.

Define an equivariant 2-form $\varpi_G\in\Omega^2_G(\g)$ by applying
the de Rham homotopy operator for the vector space $\g$ to
$\exp^*\eta_G\in\Omega^3_G(\g)$. The form degree
$0$ part of $\varpi_G$ is minus the identity map
$\g\to\g,\,Y\mapsto Y$, so that $\varpi_G(X)=\varpi-Y\cdot X$
at $Y\in\g$.
(See \cite[Section 6.2]{al:no} or \cite[Section 2]{ve:ce}
for an explicit formula for the form degree 2 part $\varpi$.)
We now define $\sigma_G=\sigma+\Psi\in\Omega^2_G(V)$ by
$$ \sigma_G=\Phi_1^*\varpi_G-\varpi_G^1-\varpi_G^2+\frac{1}{2} (\exp^*\theta^L)^1\cdot
(\exp^*\theta^R)^2.$$
Notice that $\Psi=\Phi_1-\Phi_0$ as desired, and ${\rm d}_G\sigma_G=0$ by
\eqref{eq:Mult}. Since $P_0$ vanishes at the origin, it follows
(choosing $V$ smaller if necessary)
that $1+\sigma^\flat\circ P_0^\sharp$ is invertible over $V$.
Hence \eqref{eq:change} defines a new $G$-equivariant
Poisson structure $P_1$ on $V$, with moment map
$\Phi_1(X,Y)=\log(e^Xe^Y)$.

We now scale the Poisson structure on $V$, by setting
$P_t:=t\,m_t^* P_1$ for $0<t\le 1$, with moment map $\Phi_t$.
Since $t\,m_t^* P_0=P_0$, we see that
$$ P_t^\sharp=P_0^\sharp (1+\sigma^\flat_t\circ P_0^\sharp)^{-1},$$
where $\sigma_t=t^{-1}\,m_t^*\sigma$. Since $\sigma$ is a 2-form,
$\lim_{t\to 0}\sigma_t=0$, and therefore $P_0=\lim_{t\to 0}P_t$.  Let
$\alpha_t=\frac{1}{t^2} m_t^*\alpha_1\in\Omega^1(V)$ be the family of
1-forms obtained by applying the de Rham homotopy operator to the
closed 2-forms $\frac{d \sigma_t}{d t}$. Since $P_t=t m_t^* P_1$, the Moser
vector field $v_t=P_t^\sharp(\alpha_t)$ scales according to
$v_t=\frac{1}{t} m_t^*v_1$. Note that $\alpha_t$ vanishes at the origin,
as does any form in the image of the homotopy operator.

To write $v_1$ in the form \eqref{eq:v}, define
$\g$-valued functions $A,B$ on $V$, vanishing at the origin,
by
$$ (1+\sigma_1^\flat\circ P_0^\sharp)^{-1}\alpha_1=A\cdot {\rm d} X+ B\cdot {\rm d} Y.$$
Then $P_0^\sharp(A \cdot {\rm d} X+ B \cdot{\rm d} Y)=v_1$, showing that $v_1$ is given by
\eqref{eq:v}.

\subsection{A solution to the Kashiwara-Vergne problem}
We now show that the vector field $v=v_1$ provides a solution to the
Kashiwara-Vergne conjecture in the form \eqref{eq:1'}.
The first condition $(\frac{\p}{\p t}+v_t)\Phi_t=0$
is automatic since $v_t$ is a Moser vector field.
Since the modular vector field for the
Kirillov Poisson structure of a unimodular Lie algebra vanishes, the second
condition $(\frac{\p}{\p t}+v_t)\kappa_t\Gamma=0$ holds if and only if
$v_t$ intertwines the volume forms
$\kappa_t|_{\mathcal{O}}\Lambda_{\mathcal{O}}$ on products of orbits
$\mathcal{O}=\mathcal{O}_1\times\mathcal{O}_2$, where $\Lambda_{\mathcal{O}}$ is the
symplectic volume form on $\mathcal{O}$.  By the following Proposition, the
Moser vector field $v_t$ does indeed have this property.

\begin{proposition}\label{prop:vol}
For any product of orbits $\mathcal{O}=\mathcal{O}_1\times\mathcal{O}_2\subset V$,
the volume form $\kappa_t|_{\mathcal{O}}\Lambda_\mathcal{O}$ is the Liouville
form of $\mathcal{O}$ with respect to the Poisson structure $P_t$.
\end{proposition}

Letting $\omega_{\mathcal{O}}$ denote the symplectic form on $\mathcal{O}$, this
Lemma says that the top form degree part of
$\exp(\omega_{\mathcal{O}}+\iota_{\mathcal{O}}^* \sigma_t)$ equals the top form
degree part of $\kappa_t|_{\mathcal{O}} \exp(\omega_{\mathcal{O}})$. It is
possible, but tedious, to prove this by a direct check of the
equality $\kappa_t={\det}^{1/2}(1+\sigma^\flat_t\circ P_0^\sharp)=\lambda_t$.  An
alternative route uses the theory of group-valued moment maps
\cite{al:mom}. We briefly summarize the features of this theory needed
here:

(i) Suppose $M$ is a $G$-manifold, and $\omega_G=\omega-\Phi$ an equivariant cocycle
with $\omega$ non-degenerate. Let $\Lambda$ be the Liouville form defined as
the top form degree part of $\exp\omega$.
Let $\ti{\omega}=\omega_G-\Phi^*\varpi_G$, and
$\ti{\Phi}=\exp\circ \Phi:\,M\to G$. Then $\ti{\omega}$ is non-degenerate
near $\Phi^{-1}(0)$, and $\ti{\Phi}$ satisfies the moment map condition
for group-valued moment maps, ${\rm d}_G\ti{\omega}+\ti{\Phi}^*\eta_G=0$.
Let $\ti{\Lambda}$ be the top form degree part of $\exp(\ti{\omega})$,
divided by the function $\det^{1/2}(\frac{1}{2}({\rm Ad}_\Phi+1))$. Then
$ \ti{\Lambda}=(\Phi^*J^{1/2})\ \Lambda.$

(ii) Suppose $M_1,M_2$ are two $G$-manifolds, equipped with equivariant
maps $\ti{\Phi}_j:\,M_j\to G$ and invariant 2-forms $\ti{\omega}_j$
satisfying ${\rm d}_G\ti{\omega}_j+\ti{\Phi}_j^*\eta_G=0$. Suppose for
simplicity that $\ti{\omega}_j$ are non-degenerate. Let $\ti{\Lambda}_j$
be defined as in the last paragraph. Let $M_1\times M_2$ be equipped
with the diagonal $G$-action, group-valued moment map
$\ti{\Phi}=\ti{\Phi}_1\ti{\Phi}_2$, and
2-form $\ti\omega=\ti{\omega}_1+\ti{\omega}_2+ \frac{1}{2} \ti{\Phi}_1^*\theta^L\cdot
\ti{\Phi}_2^*\theta^R.$
Then ${\rm d}_G\ti{\omega}+\ti{\Phi}^*\eta_G=0$, and the volume form $\ti{\Lambda}$
constructed from $\ti{\omega},\ti{\Phi}$ is simply the direct product
$\ti{\Lambda}_1\times\ti{\Lambda}_2$.

These facts are proved in \cite{al:du}. (The results in that paper
were stated for $\g$ is compact, but the proofs work for quadratic
$\g$, with obvious modifications.)  In our case, we start out with two
coadjoint orbits $\mathcal{O}_j$, with moment maps $\Phi_j$ the inclusion
into $\g$. Let $\Lambda_j$ be the Liouville forms.  ``Exponentiating''
as in (i), the group-valued moment maps $\ti{\Phi}_j$ are $\exp\circ
\Phi_j$, and the volume forms become $J^{1/2}\Lambda_j$. Taking the
product as in (ii), $\mathcal{O}_1\times\mathcal{O}_2$ becomes a space with
group-valued moment map $\ti{\Phi}(X,Y)\mapsto e^X\,e^Y$ and volume
form $J^{1/2}(X)J^{1/2}(Y)\Lambda$ where
$\Lambda=\Lambda_1\times\Lambda_2$.  Using (i) in reverse,
$\mathcal{O}_1\times\mathcal{O}_2$ becomes a Hamiltonian $G$-space in the usual
sense, with moment map $\Phi(X,Y)=\log(e^X e^Y)$ and equivariant
symplectic form $\omega_G+(\Phi_1,\Phi_2)^*\sigma_G$. The Liouville volume
form is $\frac{J^{1/2}(X)J^{1/2}(Y)}{J^{1/2}(Z)}\Lambda$.  This yields a
proof of Proposition \ref{prop:vol} for $t=1$, the general case
follows by rescaling the given bilinear form on $\g$.

\begin{remark}
Our solution to the Kashiwara-Vergne problem for quadratic Lie algebras
differs from Vergne's solution in \cite{ve:ce}, although the
ingredients are very similar. A common feature of both solutions is that
they are, in fact, independent of the quadratic form.
\end{remark}

%%%%%%%%%%%%%%%%%%%%%%%%%%%%%%%%%%%%%%%%%%%%%%%%%%%%%%%%%%%%
%%%  Bibliography  %%%
%%%%%%%%%%%%%%%%%%%%%%%

%
\end{document}